# On the Neuro Response Feature of Deep CNN for Remote Sensing Image

Jie Chen, Min Deng, Haifeng Li

School of Geosciences and Info Physics, Central South University, Changsha, Hunan, China
lihaifeng@csu.edu.cn

Abstract

本文拟对深度卷积网络的内在性质进行研究，希望通过分析深度卷积网络响应的规律来更好的理解深度神经网络的工作基本原理

## 1. Introduction

我国高分辨率对地观测系统重大专项已全面启动，高空间、高光谱、高时间分辨率和宽地面覆盖于一体的全球天空地一体化立体对地观测网逐步形成，将成为保障国家安全的基础性和战略性资源(李德仁等 2012，何国金等 2015)。未来 10 年全球每天获取的观测数据将超过 10PB，遥感大数据(Baraniuk 2011；李德仁，张良培等 2014；陈刚 2015)时代已然来临。Science 杂志认为遥感信息处理的现状可描述为："data-rich but analysis-poor" (Clery 2005)，也即"大数据，小知识"(李国杰 2012；李德仁，张良培等 2014)。随着小卫星星座的普及，我们可以实现每天 3 次以上的全球覆盖能力，遥感大数据与知识转化及服务能力不足的矛盾将显得越发突出，其关键理论挑战和技术瓶颈在于遥感影像的自动分析与理解依然是国内外遥感领域共同面临的开放科学问题。

遥感影像地物识别的核心问题是如何建立一个"好特征模型"来刻画地物光谱与空间影像特征。这些特征模型包括：光谱信息内容模型(Zhang et al. 2015)、纹理基元森林(Zhen et al. 2014)、规则形状相似指数(Sun et al. 2015)、形态建筑指数(Chen et al. 2014)、形态属性剖面(Xia et al. 2015)、空间方位语义特征(陈杰等 2016)、增强标记最小生成树(Akbari et al. 2016)、条件随机场(Zhang et al. 2015)、多尺度梯度方向直方图(Cheng et al. 2013)，块旋转不变模型(Zhang et al. 2014)、Gabor(Risojevi et al. 2013)、显著性(Zhang et al. 2015)。此外，为缩小机器视觉上的"语义鸿沟"，中层特征模型也在遥感影像分类与目标识别研究领域得到了快速发展。如空间金字塔匹配(Mei et al. 2015)、词袋模型(Hu et al. 2015)、主题袋模型(Bahmanyar et al. 2015)、综合光谱、纹理及 Sift 的概率潜在语义模型(Zhong et al. 2015)、狄利克雷多主题模型(Zhao et al. 2015)、稀疏表达(Cheng et al. 2015; Tao et al. 2015)。虽然人们对遥感影像自动目标识别取得不俗的进展，**但是它仍是公认的世界性开放难题**。

最近借助深度学习方法(Hinton 2006)，遥感影像自动地物识别也取得了令人印象深刻的结果(Chen et al. 2014; Hu et al. 2015; Zhang et al. 2015; Zhao et al. 2015)。深度学习*采用"端对*

端"的特征学习,通过多层处理机制揭示隐藏于数据中的非线性特征,能够从大量训练集中自动学习全局特征(*这种特征被称为"学习特征"*),是其在遥感影像自动目标识别取得成功的重要原因(Zhao et.al. 2016),也标志特征模型从**手工特征**向**学习特征**转变。更令人惊讶*的是*,在一些特定的领域和特定的数据集上,深度学习方法的表现甚至超过了人类,如 2015 年 DeepID2 项目人脸识别率准确率提升到 99.15%(He et al. 2015; Ouyang et al. 2015);微软亚洲研究院在 2015 年的 ImageNet 挑战赛中,将图像分类错误率降到了 3.5%(He et al. 2015);2016 年 Google Deep Mind 团队的 AlphaGo,在围棋对战领域,凭借强大的学习能力战胜了世界冠军李世石(Silver et al. 2016)。正是深度学习在视觉领域和智能领域超乎寻常的学习能力,引发了 Nature 和 Science 两大顶级杂志人工智能的专刊研究(Jordan 2015; LeCun et al. 2015),这一切标志着"感知和计算智能"新纪元的开启。

如何有效地训练很深的网络模型仍是未来研究的一个重要课题。一些人认为深度学习的成功在于用具有大量参数的复杂模型去拟合数据集。这个看法也是不全面的。事实上,进一步的研究表明 DeepID2+的特征有很多重要有趣的性质。例如,它最上层的神经元响应是中度稀疏的,对人脸身份和各种人脸属性具有很强的选择性,对局部遮挡有很强的鲁棒性。而 DeepID2+通过大规模学习自动拥有了这些引人注目的属性,其背后的理论分析值得未来进一步研究。虽然深度学习在实践中取得了巨大成功,通过大数据训练得到的深度模型体现出的特性,例如分类精度与深度层次的关系、滤波器对不同对象的响应情况等,都需要作深入的研究。

基于上述研究,本文拟对深度卷积网络的内在性质进行研究,希望通过分析深度卷积网络响应的规律来更好的理解深度神经网络的工作基本原理。

## 2. Experiment and parameter setting

本研究基于美国加利福利亚大学的 UC Merced 数据集(简称为 UCM)进行,它由 21 类不同土地利用类型的遥感影像数据组成,每类各包含 100 幅大小为 256×256 像元的影像块,分辨率为 0.33 米。所采用的深度学习网络以当前在计算机视觉及遥感影像处理领域使用较为普遍的 Caffe 为框架。该深度学习框架给出了模型定义、最优化设置与预训练权值,且采用模块化设计思想,便于用户根据提供的各层类型自定义模型,总体上构架清晰且计算高效,易于操作。由于本研究是通过高分辨率遥感影像的深度学习实验,揭示深度神经网络对于影像特征的逐层抽取工作机制,采用互联网上开源的两个成熟的深度学习网络 AlexNet 和 CaffeRefNet 进行实验。它们的网络主体结构是一致的,均由若干个卷积层(Conv)、池化层(Pooling)与局部响应归一化层(LRN)组成,仅是两者的 LRN 层与 Pooling 层前后次序是相反的。本文实验将对它们设计 3 层、4 层、5 层不同卷积层数的网络深度,分别记

为 AlextNet-3/4/5Conv 和 CaffeRefNet-3/4/5Conv。

具体的实验操作，是在 Ubuntu 平台上使用 Caffe 框架，主要参数设置为 stride=2、std=0.01，将 UCM 中每类影像的 80%作为训练样本，余下 20%作为测试数据，输入 227×227×3 大小的影像经 20K 次迭代训练得到稳定网络参数。

## 3. The strategies and methods

揭示深度神经网络对不同地物类型遥感影像的工作机制，需要分析某地类影像在深度神经网络中的滤波器响应规律。为此，我们分别提取出 AlextNet-3/4/5Conv 和 CaffeRefNet-3/4/5Conv 的 6 个网络结构中的最高层滤波器，并设计以下策略来分析判定某个滤波器是否对某地物类型影像具有稳定响应。

首先，基于上述 6 个网络结构从中选择测试精度较高的网络进行影像响应特征的分析。在本实验室中它们分别为 Forest、Harbor、Overpass、Sparse-Residential、Parkinglot 和 Chaparral，并根据特征接近程度分别将这 6 类划分为 6 个小类，每小类包含 3~6 张影像。

其次，只有当某个滤波器对上述某类地物中的 6 小类影像都存在响应时，才判定该滤波器与该类地物影像之间存在稳定响应，称其为稳定滤波器；若该滤波器对该类地物影像中的任意一副影像不存在响应，则将该滤波器的响应视为非稳定。在本实验中，由于实验数据量的关系，在统计每类地物影像稳定滤波器时，可能存在未达到最优化的结果，即虽然得到的某个稳定滤波器对当前地物类型实验中影像都存在响应，但是不可确定其对所有该类地物影像都存在响应。为此，实验中我们设计了稳定响应滤波测试影像集，对统计结果进行二次验证判别，尽可能提高实验统计精度。

最后，采用深度网络特征可视化工具 Deep Visualization Toolbox（DVT），直观地显示各层深度神经网络学习到的影像特征。也就是，对该 6 类影像在前述 6 种不同网络的最高卷积层中滤波器的响应，通过解卷积方式将滤波器响应特征图进行可视化。在此基础上，探究深度学习模型各层网络中滤波器卷积得到的图像特征与影像上某类地物的特征是否存在对应关系。

下面，我们将从 4 个方面对实验中得到的遥感影像深度学习特征响应规律进行描述与分析。

| 序号 | 类名 | 正确分类率 | | | | | |
| --- | --- | --- | --- | --- | --- | --- | --- |
| | | AlexNet-5Conv | AlexNet-4Conv | AlexNet-3Conv | CaffeRefNet-5Conv | CaffeRefNet-4Conv | CaffeRefNet-3Conv |
| 1 | 'chaparral' | 1 | 1 | 1 | 1 | 1 | 1 |
| 2 | 'forest' | 1 | 1 | 1 | 1 | 1 | 1 |
| 3 | 'harbor' | 1 | 0.95 | 1 | 1 | 1 | 1 |
| 4 | 'overpass' | 0.85 | 0.85 | 0.95 | 0.9 | 0.75 | 0.85 |
| 5 | 'parkinglot' | 0.6 | 1 | 0.95 | 0.75 | 1 | 1 |
| 6 | 'sparse-residential' | 0.8 | 0.85 | 0.8 | 0.85 | 0.65 | 0.95 |
| | 21类总分类精度 | 65.53% | 72.89% | 70.00% | 69.47% | 72.37% | 70.26% |

图 1 六类地物影像在 6 中网络结构中的分类精度

## 4. The flexibility of network depth

**观察一：同一结构的 DL 网络，不同地类的影像分类所需的较优网络深度不同。**

利用前述 6 种深度学习网络结构对 UCM 的 21 类影像进行分类测试，并统计它们的分类精度，结果如图 2，3。在图 2 中，对于 AlexNet 网络模型，Chaparral、Forest、Harbor、Overpass、Beach、Intensive-Residential、Freeway 等地类在网络结构 3Conv 中可以得到较高的分类精度，除 Buildings、Baseball-Diamond 在网络结构 5Conv 上分类精度最高外，其他地类的分类精度均是在 4Conv 网络结构中较高。对于 CaffeRefNet 网络，Plowland、Baseball-Diamond、Beach、Freeway、Golf-Course、Medium-Residential、Portable-House、Tennis-Court 在网络结构 4Conv 上分类精度达到峰值，也就是说，随网络深度的增加其分类精度反而会降低；Airport、Sparse-Residential 类在网络结构 3Conv 上的分类精度最优，而 Beach、Portable-House 类在 4Conv 网络上的分类精度最优。仅有 AlexNet 网络中的 Buildings 与 CAffeRefNet 网络中的 Dense-residential 的分类精度，表现为随网络深度增加而提高，其他地类如 Chaparral、Forest、Harbor、Sparse-Residential、River 等的分类精度在 AlextNet 和 CaffeRefNet 随着网络深度的变化并不明显。另外，同样深度的 AlexNet 与 CAffeRefNet 网络，同时对某些地物会表现出较好的分类精度，对有些地物的分类精度较差，如图 3 所示。其中，AlexNet-5Conv 和 CaffeRefNet-5Conv 网络均对 ID1、2、3、4、6 地物分类精度达到 80%以上，而对 ID14、16、17、18、20 地物分类精度在 50%以下。

从上述分类精度的变化中可发现：1）相同地物的最佳分类精度会在 AlexNet 或 CaffeRefNet 网络结构的不同网络深度得到；2）地物分类的精度并不随网络结构深度的增加而提高；3）AlexNet 或 CaffeRefNet 的不同网络深度，对于地物分类精度的影响总体而言并不大；4）深度学习网络的深度主要决定于地物影像特征本身。

综上所述，在不同深度学习网络结构中，UCM 的 21 类地物影像的分类精度未见随网络深度变化而呈现明显的升降规律，但网络深度的确会影响地物的分类精度，而这种影响背后的数学关系还有待研究。

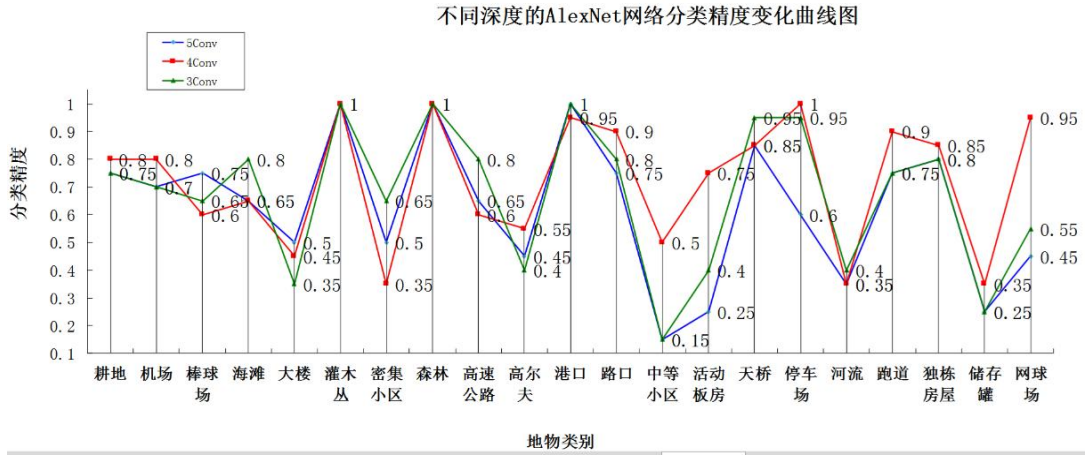

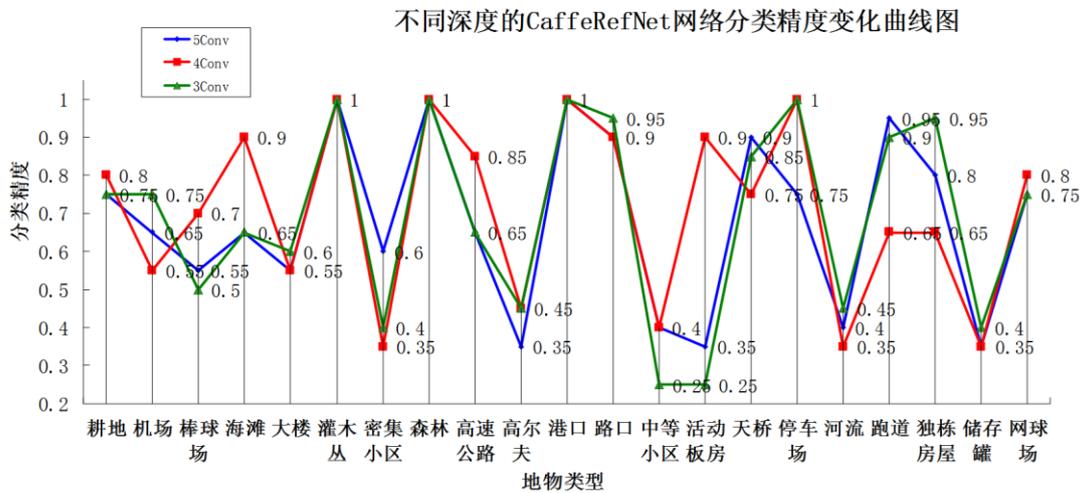

图 2 UCM 数据集在不同深度、相同结构的 DL 网络上的分类精度曲线图。将较优网络层数作为变量进行数值量化，3Conv 则赋值 3,4Conv 赋值 4,5Conv 赋值 5，则随着地类变化，各类地类较优网络的情况如图 3。

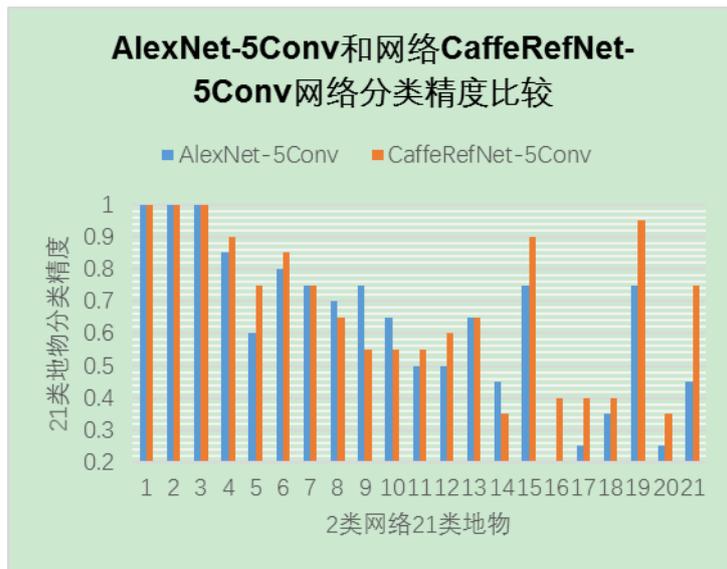

图 3 21 类地物在各层网络中的分类精度

## 5. The pattern of filter activations

**观察二：影像地物类型决定网络最高层卷积响应滤波器的组合方式**

使用 DVT 显示 6 类影像分别在前述 6 种网络中最高层滤波器的响应，并统计其中稳定响应滤波器的数目和种类 ID（最高层滤波器总数为 512，它们的种类 ID 范围取 0-511）。从数量上进行统计的结果如表 2 所示，其中的数字表示稳定滤波器的数量，而在不同深度网络中的稳定滤波器数量共同构成稳定响应的组合模式。由表中数字可知，在利用深度神经网络进行遥感影像地物分类时，最高层的一定数量的卷积响应滤波器确实对特定地类存在着稳定响应，且对于不同类地物的该类滤波器数目各不相同。以表 2 中 AextNet-4Conv 为例，各类地物在最高卷积层 Conv4 中的稳定响应滤波器数量是不同的。而从种类 ID 角度上，滤波器也表现出对于特定地物的稳定响应。以 Forest 为例（图 7），稳定响应滤波器组合由 5 种组成，对应其 ID 分别为 205、299、392、410 和 497；而 Parkinglot 类影像对应的稳定响应滤波器组合则由 7 种滤波器组成，其 ID 分别为 114、268、276、334、338、385、431；Harbor 类影像对应的滤波器组合则由 3、10、74、184、233、322、327、350、362、368、379、417 和 499 共 13 种滤波器组成。其他类影像在在网络高层卷积滤波器的响应也存在类似规律，例如在 AlexNet-3Conv 网络中，Sparse-Rsidential 最高层卷积的稳定响应由 14 种滤波器组成（图 8）。

| 类名 | 网络结构 | | | | | |
|---|---|---|---|---|---|---|
| | AlexNet-5Conv | AlexNet-4Conv | AlexNet-3Conv | CaffeRefNet-5Conv | CaffeRefNet-4Conv | CaffeRefNet-3Conv |
| chaparral | 127 | 11 | 8 | 240 | 30 | 29 |
| forest | 111 | 67 | 20 | 225 | 6 | 5 |
| harbor | 215 | 145 | 21 | 314 | 83 | 72 |
| overpass | 260 | 1 | 3 | 255 | 33 | 64 |
| parkinglot | 255 | 7 | 2 | 153 | 21 | 34 |
| sparse-residential | 140 | 17 | 14 | 340 | 121 | 127 |

图 4 6 类实验影像在 6 种网络结构中的稳定响应滤波器组合

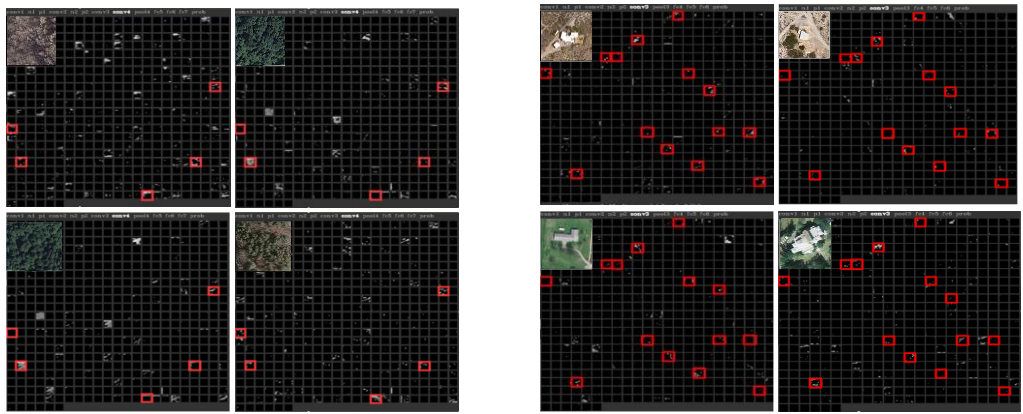

图 5 Forest 类影像在 AlexNet-4Conv 中的稳定响应的滤波器组合；图 8 Sparse-Residential 在 AlexNet-3Conv 中的稳定响应滤波器组合

((红色框位置对应14 个稳定响应滤波器的位置)(红色框位置对应 5 个稳定响应滤波器的位置)

由上述实验结果可见，相同地类影像在同一深度学习网络结构中的稳定滤波器响应具有数量与种类的组合模式，而对于不同类型地物而言，不仅其稳定滤波器的数目有差别，而且稳定滤波器的组合方式也有明显区别。

特别的我们用 t-SNE 降维可视化方法将学习特征分布投影到二维和三维空间，如图 **6** 所示，图 **7** 所示，6 类地物能够被分类器区分开来的基本原因是，因为他们的学习特征在低维空间呈现了显著的聚集现象。

这个发现有利于在用深度学习网络模型进行遥感影像分类时，获知其有效作用区域以及层间关联，这将为今后更进一步地就某个特定任务改进网络模型提供了重要的参考依据。

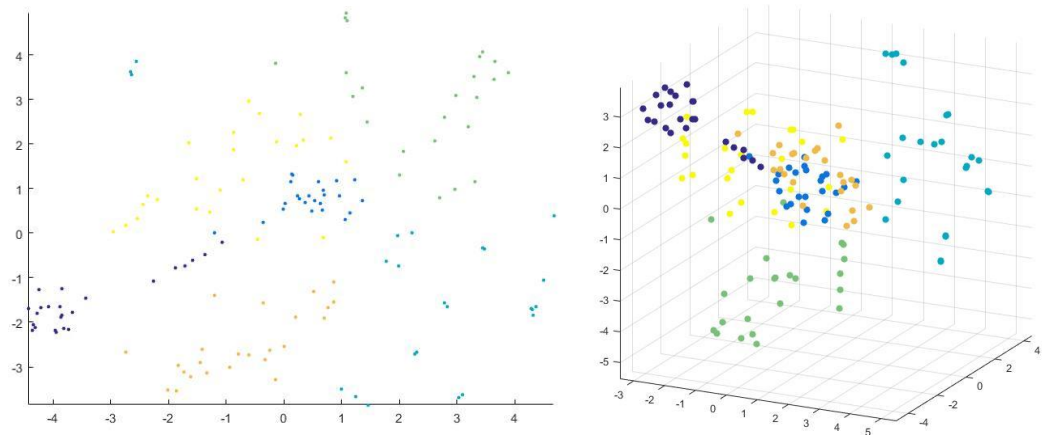

图 6 地物的特征空间分布图

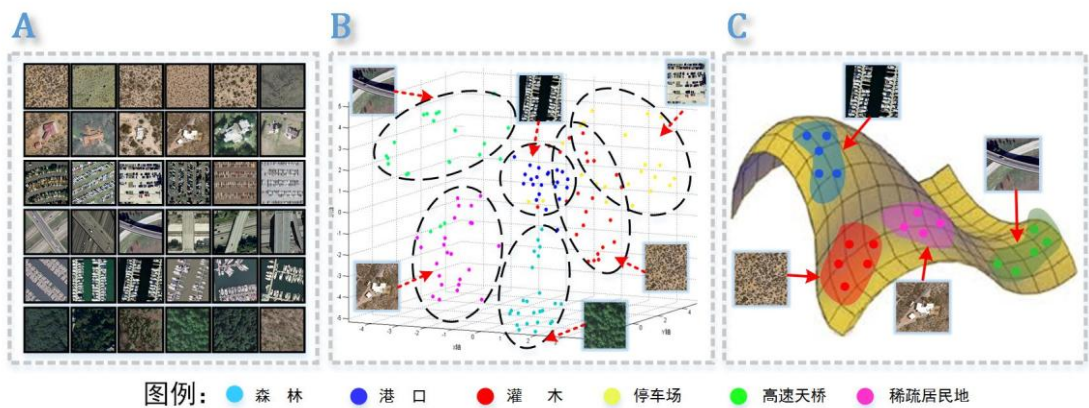

图 7 地物的特征空间分布图

## 6. Selectiveness on identities and attributes

**观察三：网络最高层卷积中存在稳定响应的滤波器，其对地物特征的响应具有选择性**

使用 DVT 显示 6 类影像分别在前述 6 种网络中最高层滤波器的响应并作解卷积处理，通过人眼判读的方式分析各高层卷积滤波器学习到的地物特征。部分实验结果如图 8 所示，

在 AlexNet-3Conv 网络的 Conv3 层中 ID 为 78 的滤波器总能对 Sparse-Residential 类的独栋房屋区域和 Harbor 类的船舶区域固定响应；ID 为 424 的滤波器则对 Harbor 类影像的港口海水区域响应良好。这种规律同样在 AlexNet-4/5Conv 网络中普遍存在，如 Conv5 中 ID 为 164 的滤波器对 Overpass 类影像的天桥区域响应稳定；Conv4 中 ID 为 417 的滤波器与 ID 为 233 号滤波器分别对 Harbor 类影像的港口海水区域和船舶区域有着固定响应；Conv4 中 ID 为 344 滤波器对 Airplane 类的飞机机体区域有良好响应等。由该部分实验结果可以发现，最高层卷积中的部分滤波对学习到的特征具有明显的选择性。

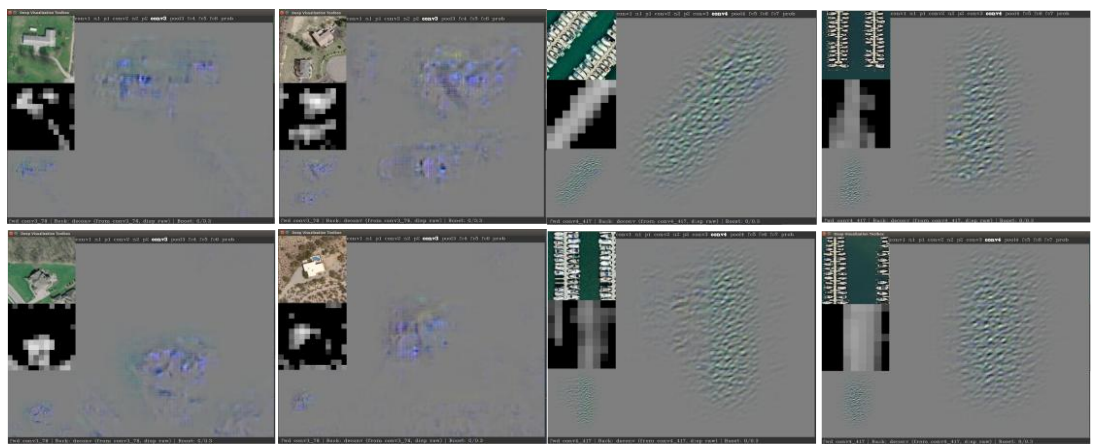

（a） （b）

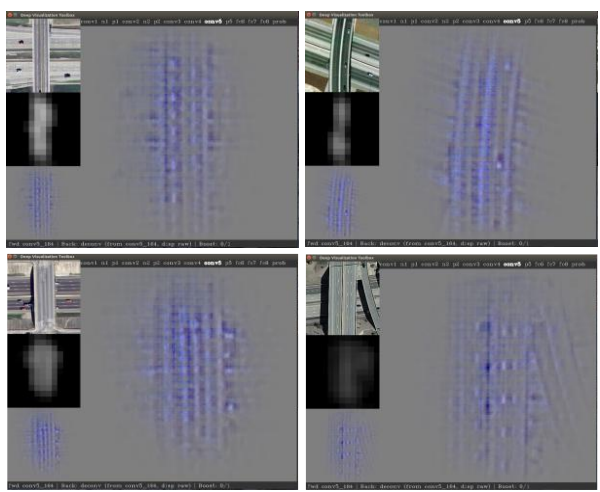

(c)

（a）AlexNet-3Conv 网中 78 号滤波器对 Sparse-Residential 影像独栋建筑物的响应
（b）AlexNet-4Conv 网中 417 号滤波器对 Harbor 影像的港口区域的响应
（c）AlexNet-5Conv 网中 164 号滤波器对 Overpass 影像的过街天桥区域的响应
图 8 网络高层卷积滤波器对地物特征的选择性响应示意图

## 7. Robustness of the deep features

**观察四：不同地类的遮挡图像在卷积神经网络模型下鲁棒性实验规律**

将 UCM 中全部 2100 幅影像作为分类测试数据，在 AlexNet 网络模型上对不同地类影

像开展深度神经网络的特征鲁棒性测试。具体地，分别对 21 类地物的 100 幅影像进行部分遮挡，并统计各类地物遮挡后各类影像的分类精度。遮挡形式分两种：比例遮挡（10%-50%）与随机遮挡（遮挡块大小依次为 13×13、26×26、38×38、51×51、77×77、102×102），如图 9 所示。

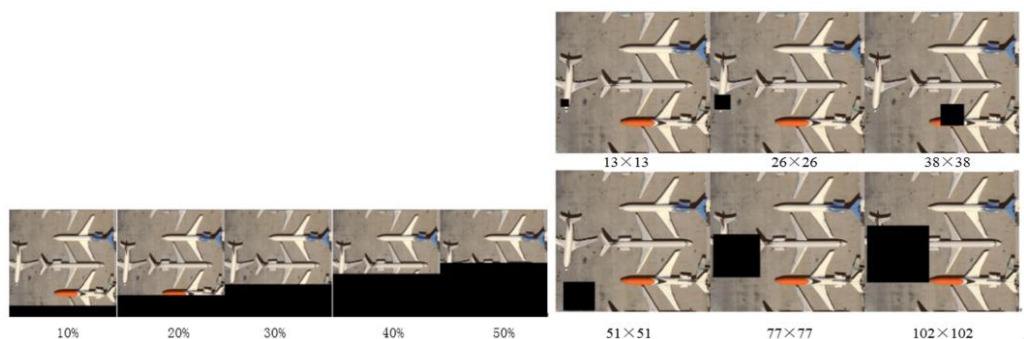

（b）机场图像的 10%-50%的百分比遮挡
（a）机场图像六种的随机块遮挡
图 9 实验采用的遮挡测试集示意图（以机场为例）

对基础测试集进行百分比遮挡和指定大小的块遮挡后，在 AlexNet-Conv3、AlexNet-Conv4、AlexNet-Conv5 网络的训练模型上，分别对两个遮挡测试集进行测试。从如图 10 所示结果中可以看出，表现了不同深度 AlexNet 网络训练模型在 UCM 遮挡测试集上随遮挡的分类精度变化情况。在遮挡率低于 15%时，三类网络的分类精度都在 50%以上，在遮挡率高于 30%时，分类精度急剧降低。三类网络中，AlexNet-5Conv 网络模型精度趋势相比于其他两类网络下降较为平缓，AlexNet-4Conv 相比于 AlexNet-3Conv 网络精度下降较为平缓。

在指定遮挡块的遮挡实验中，得到不同深度的网络在 UCM 遮挡测试集上分类总体精度随遮挡的变化，如图 11。在随机遮挡块小于 38×38 时，三类网络模型的总体精度仍能保持在 50%以上，其精度下降较为平缓，在遮挡块大于 51×51 时，分类精度便急剧下降。与百分比遮挡变化的结果类似的是，AlexNet-5Conv 网络模型在三类网络中分类总体精度的变化率最低，AlexNet-4Conv 网络次之，AlexNet-3Conv 网络变化趋势最明显。这说明三类网络中，AlexNet-5Conv 网络具有较强的鲁棒性，对此，我们推测较深的网络可能具有更强的鲁棒性。

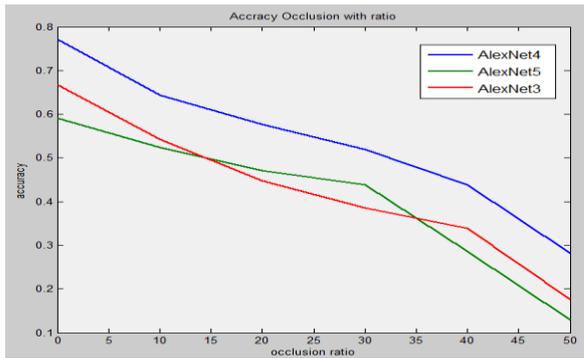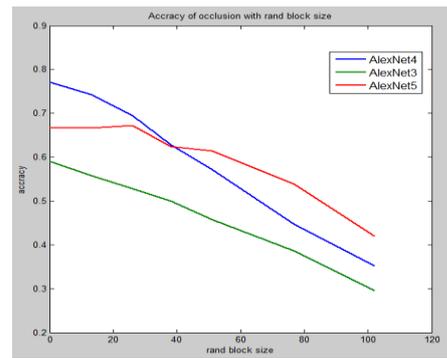

图 10 UCM 遮挡集上分类总体精度随着遮挡百分率的变化（x 轴表示遮挡部分占全图百分率，y 轴表示分类的总体精度，）

图 11 UCM 测试遮挡集上不同深度 AlexNet 网络分类总体精度随遮挡块的变化（x 轴表示遮挡块的大小，y 轴表示分类的总体精度）

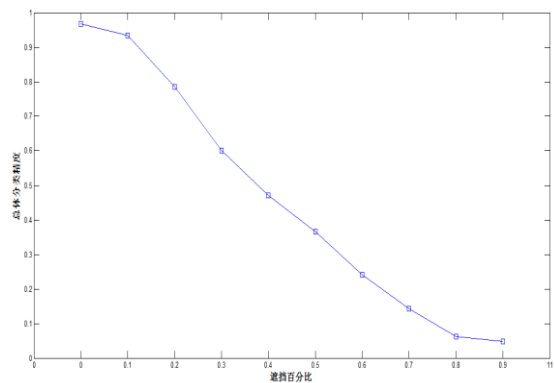

图 12 在 AlexNet 网络下的整体分类精度随遮挡比例的变化曲线图

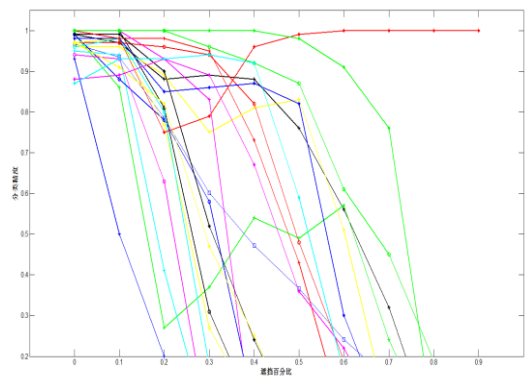

图 13 不同地类的高分影像在 AlexNet 网络中的分类精度随遮挡率的变化曲线图

图 **12** 为 2100 幅影像的整体分类精度随遮挡百分比增加的变化曲线图，图 **13** 给出了 21 类地物的单类分类精度随遮挡百分比增加的变化曲线图。由图 1 可见，随遮挡百分比的增加，AlexNet 网络模型的分类精度急剧降低。在遮挡比例超过了 0.2 以后，该网络模型的整体分类精度呈线性递减趋势。而由图 2 可知，不同地物类型在该网络模型下的分类精度随遮挡比例的变化趋势并不一致。但仍可明显看出 0.2 的遮挡比例是一个阈值点，绝大多数地类的分类精度在遮挡比例超过 0.2 之后骤然降低，而在 0.2 的遮挡比例之前，分类精度相对降低幅度并不很大。这种特性在一定程度上证明了卷积神经网络提取的高层图像特征具有一定程度

的抗遮挡能力。另外，仔细观察上图可知，在 UCM 数据集中有几类地物场景（如河流、跑道、港口等）的抗遮挡能力明显较强，当遮挡比例超过了 0.4-0.5 时其分类精度仍然可以维持在一个较高的水平。主要因为这几种地物的内部结构相对均一，局部遮挡对整体影像分类的影响相对较小。

## 8. Conclusion

通过深入分析深度神经网络这个黑盒子内部，我们发现深度神经网络背后蕴含这非常有趣的规律：深度伸缩性，特征选择性，神经元组合模式稳定性和识别鲁棒性。这些特性可以为研究人员更好的使用深度神经网络提供参考。